\documentclass[11pt,a4paper]{amsart}
\usepackage{amssymb, amstext, amscd, amsmath}
\usepackage{url}
\usepackage{amsfonts}
\usepackage{lscape}
\usepackage{amsthm}
\usepackage{amsfonts}
\usepackage{array}
\usepackage{latexsym}
\usepackage{mathrsfs}
\usepackage{textcomp}
\usepackage{verbatim}
\usepackage{tikz}
\usepackage{blkarray}
\usepackage{thmtools, thm-restate}

\newtheorem{thm}{Theorem}[section]

\newtheorem{claim}[thm]{Claim}

\newtheorem{cor}[thm]{Corollary}
\newtheorem{lem}[thm]{Lemma}

\theoremstyle{definition}
\newtheorem{example}{Example}[section]

\numberwithin{equation}{section}

\newcommand{\R}{{\mathbb{R}}}

\newcommand{\bZ}{{\mathbb{Z}}}

  \newcommand{\I}{{\mathcal{I}}}

  \newcommand{\M}{{\mathcal{M}}}

  \newcommand{\calR}{{\mathcal{R}}}

  \newcommand{\X}{{\mathcal{X}}}

\newcommand{\rank}{\operatorname{ rank}}

\newcommand{\val}{\mathrm{val}}

\usetikzlibrary{patterns}
\tikzset{roundnode/.style={circle,draw=black!50,fill=black!20,inner sep=1.2pt}}
\tikzset{every loop/.style={}}
\tikzset{polarshift/.style args={#1/#2}{xshift=#1*cos(#2),yshift=#1*sin(#2)}}
\begin{document}
\pagestyle{plain}

\title{A sufficient connectivity condition for rigidity and global rigidity of linearly constrained frameworks in $\R^2$ }

\author[Hakan Guler]{Hakan Guler}
\address{Department of Mathematics, Faculty of Sciences, Kastamonu University, Kastamonu, 37150, T\"urkiye}
\email{hakanguler19@gmail.com}
\email{hguler@kastamonu.edu.tr}

\date{\today}

\begin{abstract}
We study the bar-and-joint frameworks in $\R^2$ such that some vertices are constrained to lie on some lines.
The generic rigidity of such frameworks is characterised by Streinu and Theran \cite{ST}.
Katoh and Tanigawa \cite{KatTan} remarked that the corresponding matroid and its rank function can be
characterised by using a submodular function. In this paper, we will transfer this characterisation of the rank function to the form of the
value of a ``1-thin cover" and obtain a sufficient connectivity condition for rigidity and global rigidity of these frameworks
analogous to the results of Lov\'asz and Yemini \cite{LY}.
\end{abstract}


\maketitle

\section{Introduction}
A (2-dimensional) {\em bar-and-joint framework} is a pair $(G,p)$ where $G=(V,E)$ is a simple graph and $p:V\rightarrow \R^2$ is a map
which is called the {\em realisation} map of the framework $(G,p)$. The framework $(G,p)$ is called {\em rigid}
if every continuous motion which fixes the edge lengths is a congruence of $\R^2$. Determining whether a framework $(G,p)$ is rigid
is NP-hard \cite{A}. The problem is easier to deal with when $p$ is {\em generic}, that is the multiset of coordinates of $p$ is algebraically independent over $\mathbb{Q}$.
The rigidity of such frameworks only depends on the underlying graph $G$ \cite{AR}.

\begin{figure}[h]
\begin{center}
\begin{tikzpicture}[scale=1.3,]
\node[roundnode] at (1,0) (a) [label=below:$u$]{}
	edge[in=0,out=320,loop] node[pos=.5,below]{$l$} ()
	edge[in=40,out=0,loop] node[pos=.5,above]{$s$} ();
\node[roundnode] at (-1,0) (b) [label=below:$v$]{}
	edge[] (a) ;
\node[roundnode] at (0.5,1) (c) [label=left:$w$]{}
	edge[] (a)
	edge[](b)
	edge[in=70,out=110,loop] node[pos=.6,right]{$t$}();
\node[roundnode] at (-0.5,1) (c) [label=left:$x$]{}
	edge[in=70,out=110,loop] node[pos=.6,right]{$z$}();
\node[] at (0,0) (G) [label=below:$G$] {};
\node[] at (2,1.4) () [label=right:$p(u)\text{=}(u_1\text{,}u_2)$]{};
\node[] at (2,1) () [label=right:$p(v)\text{=}(v_1\text{,}v_2)$]{};
\node[] at (2,0.6) () [label=right:$p(w)\text{=}(w_1\text{,}w_2)$]{};
\node[] at (2,0.2) () [label=right:$p(x)\text{=}(x_1\text{,}x_2)$]{};

\node[] at (4.5,1.4) () [label=right:$q(l)\text{=}(l_1\text{,}l_2)$]{};
\node[] at (4.5,1) () [label=right:$q(s)\text{=}(s_1\text{,}s_2)$]{};
\node[] at (4.5,0.6) () [label=right:$q(t)\text{=}(t_1\text{,}t_2)$]{};
\node[] at (4.5,0.2) () [label=right:$q(z)\text{=}(z_1\text{,}z_2)$]{};
\end{tikzpicture}
\end{center}
\caption{The looped simple graph $G$ and the maps $p$ and $q$ in Example \ref{exmp:linconmatrix}.}
\label{fig:lincon_RigidityMatrix}
\end{figure}
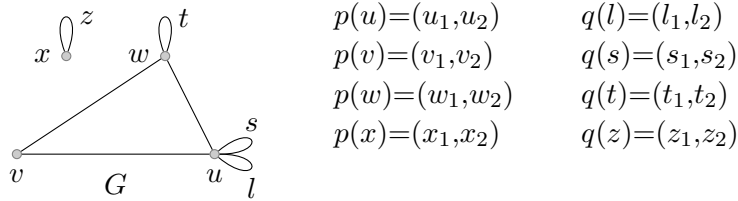

The {\em rigidity matrix} $R(G,p)$ of $(G,p)$ is a $|E|\times 2|V|$ matrix for which each edge in $E$ corresponds to a row and each vertex in $V$ corresponds to two columns
such that in the row indexed by an edge $e=uv$ we have $p(u)-p(v)$, respectively $p(v)-p(u)$ in the two columns indexed by $u$, respectively $v$ and zeros elsewhere.
We can construct a matroid $\calR(G,p)$ on $E$
from $R(G,p)$ by defining a set $F\subseteq E$ to be independent in $\calR(G,p)$ if the set of rows of $R(G,p)$ corresponding to the edges in $F$ is linearly independent.
A framework $(G,p)$ is {\em infinitesimally rigid} if $\rank R(G,p)=2|V|-3$ or equivalently if $E$ has rank $2|V|-3$ in $\calR(G,p)$.
Since the rigidity of $(G,p)$ depends only on the underlying graph $G$ when $p$ is generic, we can talk about the rigidity of a graph $G$.
A graph $G$ is {\em rigid} as a bar-and-joint framework if there exists a generic $p$ for which $(G,p)$ is infinitesimally rigid.
The graphs which are rigid as a bar-and-joint framework in $\R^2$ are characterised by Pollaczek-Geiringer \cite{PG} and rediscovered by Laman \cite{L}.
Using the fact that generic frameworks give rise to the same matroid, we can define the generic {\em rigidity matroid} $\calR_2(G)$ of $G$
by setting $\calR_2(G)=\calR(G,p)$ for some generic $p$. A characterisation of the rank function $r_2$ of $\calR_2$ is given by Lov\'asz and Yemini \cite{LY}
which is given below. The technical definitions in the theorem below will be given in Section \ref{sec:main}.
\begin{thm}\label{thm:LY}\cite{LY}
Let $G=(V,E)$ be a simple graph. Then
$$r_2(G)=\min \sum_{i=1}^k(2|X_i|-3)$$
where the minimum is taken over all 1-thin covers of $G$.
\end{thm}
With the help of Theorem \ref{thm:LY}, they also show the following.
\begin{thm}\label{thm:LY-6connected}\cite{LY}
Every 6-connected graph is rigid as a bar-and-joint framework in $\R^2$.
\end{thm}
In this paper we will consider the frameworks that have some extra constraints. These constraints will force some vertices to lie on
some lines. A {\em looped simple graph} $G=(V,E,L)$ is a graph such that $(V,E)$ is a simple graph and $L$ is a set of
loops. Note that there are no multiple edges in $E$, but $L$ is allowed to have multiple loops at a single vertex. A 2-dimensional {\em linearly constrained framework}
is a triple $(G,p,q)$ where $G=(V,E,L)$ is a looped simple graph, $p:V\rightarrow \R^2$ and $q:L\rightarrow \R^2$ are maps.
We interpret this definition as follows: $(G',p)$ where $G'=(V,E)$ (i.e.,\ $G'$ is the underlying simple graph of $G$)
is a usual 2-dimensional (bar-and-joint) framework and $q(l)$, where $l$ is a loop at $v$, constraints $v$ to lie
on the line which contains the point $p(v)$ and has $q(l)$ as its normal vector.
The {\em rigidity matrix} $R(G,p,q)$ of $(G,p,q)$ is a $(|E|+|L|)\times 2|V|$ matrix which is obtained from $R(G',p)$ by adding $|L|$ new rows
such that the row corresponding to a loop $l$ at $v$ has $q(l)$ in the two columns corresponding to $v$ and zeros elsewhere, see the example below.
\begin{example}\label{exmp:linconmatrix}
Let $G=(V,E,L)$ be the looped simple graph drawn in Figure \ref{fig:lincon_RigidityMatrix} where $V=\{u,v,w,x\}$, $E=\{uv,uw,vw\}$, $L=\{l,s,t,z\}$
and let the maps $p$ and $q$ be defined as in the figure. Then the rigidity matrix $R(G,p,q)$ is the matrix below.
\[
\begin{blockarray}{ccccccccc}
&u&u&v&v&w&w&x&x\\
\begin{block}{c[cccccccc]}
  uv & u_1-v_1 & u_2-v_2 & v_1-u_1 & v_2-u_2 & 0& 0&0 &0 \\
  uw & u_1-w_1 & u_2-w_2 & 0 & 0 & w_1-u_1 & w_2-u_2&0 &0 \\
  vw & 0 & 0 & v_1-w_1 & v_2-w_2 & w_1-v_1 &w_2-v_2&0 &0 \\
  l & l_1 & l_2 & 0 & 0 & 0 &0&0 &0 \\
  s & s_1 & s_2 & 0 & 0 & 0 &0&0 &0 \\
  t & 0 & 0 & 0 & 0 & t_1 &t_2&0 &0\\
  z & 0 & 0 & 0 & 0 & 0 &0&z_1 &z_2\\
\end{block}
\end{blockarray}
 \]
The first three rows of the matrix $R(G,p,q)$ above correspond to the rigidity matrix $R(G',p)$ where $G'=(V,E)$.
\end{example}

A linearly constrained framework in $\R^2$ is {\em infinitesimally rigid} if $\rank R(G,p,q)=2|V|$. A looped simple graph
is {\em rigid} as a linearly constrained framework in $\R^2$ if $\rank R(G,p,q)=2|V|$ for some $(p,q)$.
Similar to the bar-and-joint case we can construct a generic {\em linearly constrained rigidity matroid}
$\calR^{lc}_2(G)$ of a looped simple graph $G=(V,E,L)$ from a linearly constrained rigidity matrix $R(G,p,q)$ for which $(p,q)$ is generic.
For more information on the linearly constrained frameworks as well as some higher dimensional results
we refer the interested readers to \cite{CGJN,JNT}.

Streinu and Theran characterised the rigid looped simple graphs in $\R^2$ \cite{ST}. Before giving their result we first need to give some definitions.
Let $G=(V,E,L)$ be a looped simple graph. For a set $T\subseteq E\cup L$, we use $V(T)$
to denote the set of vertices incident with the edges and loops in $T$.
We say that $G$ is $(2,0,3)$-{\em graded-sparse} if
\begin{itemize}
\item[(i)] $|T|\leq 2|V(T)|-3$ for all $\emptyset \subsetneq T\subseteq E$ and
\item[(ii)] $|T|\leq 2|V(T)|$ for all $T\subseteq E\cup L$.
\end{itemize}
A $(2,0,3)$-graded-sparse graph is called $(2,0,3)$-{\em tight} if $|E|+|L|=2|V|$.
\begin{thm}\label{thm:ST}\cite{ST}
Let $G=(V,E,L)$ be a looped simple graph. Then $G$ is rigid (as a linearly constrained framework) in $\R^2$ if and only if
it contains a spanning $(2,0,3)$-tight subgraph.
\end{thm}

Let $r^{lc}_2$ denote the rank function of $\calR^{lc}_2$.
Katoh and Tanigawa extended Theorem \ref{thm:ST} to a non-generic setting \cite{KaTa} and pointed out that $\calR^{lc}_2$ and so $r^{lc}_2$ can be characterised by using a
submodular function \cite{KatTan}. We will transfer this characterisation of $r^{lc}_2$ to the form of the value of a 1-thin cover analogous to Theorem \ref{thm:LY}.
Since every simple graph $(V,E)$ can be considered as a looped simple graph $(V,E,L)$ with $L=\emptyset$, this
will also be an extension of Theorem \ref{thm:LY}. We will then use this extension to obtain a sufficient
connectivity condition for the rigidity and global rigidity of a looped simple graph as a linearly constrained framework in $\R^2$ and the corresponding
rigidity result will be analogous to Theorem \ref{thm:LY-6connected}.


\section{The Rank Function}\label{sec:main}

Let $X$ be a set and $g:2^X\rightarrow \bZ$ be a function. We say that $g$ is {\em submodular} if
$$g(Y)+g(Z)\geq g(Y\cup Z)+g(Y\cap Z)$$
holds for all $Y,Z\subseteq X$. 

Given a looped simple graph $G=(V,E,L)$ let us define a function $f:2^{E\cup L}\rightarrow \bZ$ by
$$
f(T):=\begin{cases}
2|V(T)|-3,& T\subseteq E\\
2|V(T)|,& T\cap L\neq \emptyset
\end{cases}
$$
where $V(T)$ is the set of all vertices incident with the edges or loops in $T$.
\begin{lem}\label{lem:submodular}
The function $f$ is submodular.
\end{lem}
\begin{proof}
Let $T_1,T_2\subseteq E\cup L$, and $a,b\in\{0,3\}$ such that we have $f(T_1)=2|V(T_1)|-a$ and $f(T_2)=2|V(T_2)|-b$. By symmetry we may assume $a\leq b$.
Then we have $f(T_1\cup T_2)=2|V(T_1\cup T_2)|-a$ since if $T_1$ or $T_2$ contains some loops then so does $T_1\cup T_2$;
and $f(T_1\cap T_2)\leq 2|V(T_1\cap T_2)|-b$ (since $T_1\cap T_2$ does not necessarily contain loops when $T_1$ and $T_2$ do so). This gives
\begin{align*}
f(T_1)+f(T_2)&=2|V(T_1)|-a + 2|V(T_2)|-b\\
			 &=2|V(T_1)\cup V(T_2)|-a + 2|V(T_1)\cap V(T_2)|-b\\
			 &\geq 2|V(T_1\cup T_2)|-a + 2|V(T_1\cap T_2)|-b\\
			 &\geq f(T_1\cup T_2)+f(T_1\cap T_2)
\end{align*}
where the first inequality follows from $V(T_1\cap T_2)\subseteq (V(T_1)\cap V(T_2))$.
\end{proof}

Katoh and Tanigawa \cite{KatTan} pointed out that the following result can be obtained by combining Lemma \ref{lem:submodular} and Edmonds \cite{Edm}.
\begin{thm}\label{thm:edmonds}\cite{KatTan}
Let $G=(V,E,L)$ be a looped simple graph and $f$ be the function defined above. Put
$$\I_f:=\{T\subseteq E\cup L: |I|\leq f(I) \text{ for all } I\subseteq T \text{ with }I\neq \emptyset\}$$
Then $\M(E\cup L, \I_f)$ is a matroid with rank function $\hat f:2^{E\cup L}\rightarrow \bZ$ given by
$$\hat f(T):=\min\big\{|T'|+\sum_{i=0}^kf(T_i):T'\subseteq T \text{ and } \{T_0,T_1,\ldots,T_k\} \text{ is a partition of }T\setminus T'\big\}.$$
\end{thm}

By Theorems \ref{thm:ST} and \ref{thm:edmonds} we have the following.
\begin{cor}\label{thm:isomorphic}
Let $G=(V,E,L)$ be a looped simple graph. Then $\M(E\cup L, \I_f)$ is isomorphic to $\calR^{lc}_2(G)$ and so
$$r^{lc}_2(T)=\hat f(T)=\min\big\{|T'|+\sum_{i=0}^kf(T_i):T'\subseteq T \text{ and } \{T_0,T_1,\ldots,T_k\} \text{ is a partition of }T\setminus T'\big\}.$$
\end{cor}
\begin{proof}
The statement follows by comparing the independent sets in each matroid, that is $(2,0,3)$-graded-sparsity is equivalent to the conditions
in the definition of $\I_f$. One may think that the empty set seems problematic in this comparison as $f(\emptyset)=-3$.
But the condition $I\neq \emptyset$ in the definition of $\I_f$ clears this issue.
\footnote{Another way of dealing with this is to set $f(\emptyset)=0$. Then $f$ will be intersecting submodular instead of submodular.
We can still use Edmonds' result \cite{Edm} on this new $f$ to obtain Theorem \ref{thm:edmonds}.
See, for example \cite[Theorem~13.4.2]{AF} for the statement of Edmonds' theorem specifically for the intersecting submodular functions.}
\end{proof}

From now on, in a looped simple graph $G=(V,E,L)$, we will sometimes need to distinguish the members of $E$ and the members of $L$.
To this end, we will call a member of $E$ a {\em (simple) edge} and a member of $L$ a {\em loop}.

Let $G=(V,E,L)$ be a looped simple graph, $T'\subseteq T\subseteq E\cup L$, and $\{T_0,T_1,\ldots,T_k\}$
be a partition of $T\setminus T'$. Consider a simple edge $e=xy\in T$. If $e\in T'$ then we can set $T_{k+1}=\{e\}$ and write

$$|T'|+\sum_{i=0}^kf(T_i)=|T'-e|+\sum_{i=0}^{k+1}f(T_i)$$
as $f(T_{k+1})=f(\{e\})=2|\{x,y\}|-3=1$. Therefore since $\{T_0,T_1,\ldots,T_k,T_{k+1}\}$ is a partiton of $T\setminus (T'-e)$,
we have the following result by applying the same operation for all simple edges in $T'\cap E$ and Corollary \ref{thm:isomorphic}.
\begin{lem}\label{lem:isomorphic_onlyloop}
Let $G=(V,E,L)$ be a looped simple graph. Then
$$r^{lc}_2(T)=\min\big\{|T'|+\sum_{i=0}^kf(T_i):T'\subseteq T\cap L \text{ and } \{T_0,T_1,\ldots,T_k\} \text{ is a partition of }T\setminus T'\big\}.$$
\end{lem}


The rank function $r^{lc}_2$ has no condition on the partition
$\{T_0,T_1,\ldots,T_k\}$ of $T\setminus T'$, in particular $|V(T_i)\cap V(T_j)|$ can be arbitrarily large and every $T_i$ is allowed to contain
some loops. We will show that the minimum value can be obtained from the partitions which satisfy $|V(T_i)\cap V(T_j)|\leq 1$ and have only one member,
say $T_0$, that may contain some loops. We need the following lemma for this.
\begin{lem}\label{lem:union}
Let $G=(V,E,L)$ be a looped simple graph and $T_1,T_2\subseteq E\cup L$.
\begin{itemize}
\item[](i) If $|V(T_1)\cap V(T_2)|\geq 2$, then $f(T_1\cup T_2)\leq f(T_1)+f(T_2)$.
\item[](ii) If $T_1\cap L$ and $T_2\cap L$ are both nonempty, then $f(T_1\cup T_2)\leq f(T_1)+f(T_2)$.
\end{itemize}
\end{lem}
\begin{proof}
We shall prove (i) and (ii) simultaneously. 
As in the proof of Lemma \ref{lem:submodular} let $a,b\in\{0,3\}$ such that we have $f(T_1)=2|V(T_1)|-a$ and $f(T_2)=2|V(T_2)|-b$.
By symmetry we may assume $a\leq b$. Then we have $f(T_1\cup T_2)=2|V(T_1\cup T_2)|-a$. Note that (i) corresponds to $b=0$ or $b=3$ and (ii)
corresponds to $b=0$. We can write
\begin{align*}
f(T_1)+f(T_2)&=2|V(T_1)|-a + 2|V(T_2)|-b\\
			 &=2|V(T_1)\cup V(T_2)|-a + 2|V(T_1)\cap V(T_2)|-b\\
			 &= 2|V(T_1\cup T_2)|-a + 2|V(T_1)\cap V(T_2)|-b\\
			 &\geq f(T_1\cup T_2)
\end{align*}
where the inequality follows from the fact that $|V(T_1)\cap V(T_2)|\geq 2$ if $b=3$, and it trivially holds when $b=0$.
\end{proof}
We can now give the rank function $r^{lc}_2$ as follows.
\begin{thm}\label{thm:better_rank_edge}
Let $G=(V,E,L)$ be a looped simple graph and $f$ be the function defined above. Then
$$r^{lc}_2(T)=\min\big\{|L'|+\sum_{i=0}^kf(T_i):L'\subseteq T\cap L \text{ and } \{T_0,T_1,\ldots,T_k\} \text{ is a partition of }T\setminus L'\big\}$$
where $T_1,T_2,\ldots,T_k$ have no loops and $|V(T_i)\cap V(T_j)|\leq 1$ for distinct $i,j$.
\end{thm}
\begin{proof}
The condition $L'\subseteq T\cap L$ follows from Lemma \ref{lem:isomorphic_onlyloop} ($T'$ is relabelled as $L'$) and the
conditions $T_1,T_2,\ldots,T_k$ have no loops and $|V(T_i)\cap V(T_j)|\leq 1$ for distinct $i,j$ follow from Lemma \ref{lem:union} and by relabelling if necessary.
(By relabelling only the set $T_0$ is allowed to contain loops.)
\end{proof}

The rank function in Theorem \ref{thm:better_rank_edge} deals with both edges or loops (partitions) and vertices (condition on $|V(T_i)\cap V(T_j)|$).
We can replace the partition of the edges or loops by some family of subsets of $V$ so that the rank function deals only with some subsets of the vertices. In order to do this we
need to give some definitions first.
Let $G=(V,E,L)$ be a looped simple graph. A family $\X=\{X_0,X_1,X_2,\ldots, X_k\}$ of subsets of $V$ satisfying $|X_i|\geq 2$, for all $1\leq i\leq k$,
is said to be a {\em cover} of $G$ if every edge $e\in E$ and every loop $l\in L$ is induced by at least one member in $\X$.
A cover $\X=\{X_0,X_1,X_2,\ldots, X_k\}$ of $G-L'$ where $L'\subseteq L$ is said to be {\em admissible} if every loop in $L\setminus L'$ is induced by the vertices in $X_0$
and every loop in $L'$ is induced by the vertices in $V\setminus X_0$.
The set $X_0$, which is allowed to be empty, is called the {\em looped member} of $\X$. We say $\X$ is {\em $t$-thin} (for some $t$) if
$|X_i\cap X_j|\leq t$ for all $0\leq i< j\leq k$. The {\em value} of the admissible cover $\X$ of $G-L'$ which is denoted by $\val(\X)$ is defined as
$$\val(\X):=|L'|+2|X_0|+\sum_{i=1}^k(2|X_i|-3).$$ For a set $X\subset V$, we use $E_G(X)$ and $L_G(X)$ to denote the set of simple edges and respectively loops
in $G$ which are induced by the vertices in $X$.

\begin{thm}\label{thm:rank_vertex}
Let $G=(V,E,L)$ be a looped simple graph. Then
$$r^{lc}_2(G)=\min\{\val(\X): L'\subseteq L, \X \text{ is an admissible 1-thin cover of } G-L'\}.$$
\end{thm}
\begin{proof}
Let $\X=\{X_0,X_1\ldots,X_k\}$ be an admissible 1-thin cover of $G-L'$ for some $L'\subseteq L$. Consider the sets $T_0=E_G(X_0)\cup L_G(X_0)$,
and $T_i=E_G(X_i)$, $1\leq i\leq k$. Note that $V(T_j)\subseteq X_j$ for all $0\leq j\leq k$. Possibly some $T_i$ are empty.
Since the empty set is not allowed in a partition, we discard all such $T_i$, and by relabelling if necessary, we may assume
$\{T_0,T_1,\ldots,T_m\}$ is a partition of $(E\cup L)\setminus L'$ for some $m\leq k$, and the sets $T_i$, $0\leq i\leq m$, satisfy the
conditions in the statement of Theorem \ref{thm:better_rank_edge}.
Now the facts that $f(T_0)\leq 2|X_0|$ (we do not know whether $T_0$ contains a loop or not), $f(T_i)\leq 2|X_i|-3$ for all $1\leq i\leq m$,
$2|X_i|-3\geq 1$, for all $1\leq i\leq k$ and Theorem \ref{thm:better_rank_edge} imply that
$\displaystyle{r^{lc}_2(G)\leq |L'|+\sum_{i=0}^mf(T_i)\leq |L'|+2|X_0|+\sum_{i=1}^k(2|X_i|-3)=\val(\X)}$.

We need to show that there exists an admissible 1-thin cover whose value hits the rank to finish the proof.
To do this let $L'\subseteq L$ and $T_0,T_1,\ldots,T_k$ be a partition of $(E\cup L)\setminus L'$ satisfying
the conditions in the statement of Theorem \ref{thm:better_rank_edge}
from which $r^{lc}_2(G)$ can be obtained i.e.,\ $r^{lc}_2(G)=|L'|+\sum_{i=0}^kf(T_i)$. We split the proof into two cases.\\
\textbf{Case 1.} $T_0\cap L\neq \emptyset$.\\
Let $\X=\{X_0,X_1,\ldots,X_k\}$ such that $X_i=V(T_i)$ for all $0\leq i\leq k$. We claim that every loop in $L'$ is induced by some vertex in $V\setminus X_0$.
To see this suppose the contrary and let $l\in L'$ be a loop at a vertex in $X_0$. Then $T_0\cup\{l\},T_1,\ldots,T_k$ is a partition of $(E\cup L)\setminus (L'\setminus\{l\})$.
Thus by Theorem \ref{thm:better_rank_edge} and the fact that $f(T_0\cup\{l\})=f(T_0)$ we have 
$$r^{lc}_2(G)\leq |L'\setminus\{l\}|+f(T_0\cup\{l\})+\sum_{i=1}^kf(T_i)=|L'|+\sum_{i=0}^kf(T_i)-1=r^{lc}_2(G)-1,$$
a contradiction.
Hence $\X$ is an admissible 1-thin cover of $G-L'$.
We have $f(T_0)=2|X_0|$ and $f(T_i)=2|X_i|-3$ for all $1\leq i\leq k$ and so

$$r^{lc}_2(G)=|L'|+\sum_{i=0}^kf(T_i)=|L'|+2|X_0|+\sum_{i=1}^k(2|X_i|-3)=\val(\X).$$
\textbf{Case 2.} $T_0\subseteq E$.\\
Let $\X=\{X_0,X_1,X_2,\ldots,X_k,X_{k+1}\}$ where $X_0=\emptyset$, $X_i=V(T_i)$ for all $1\leq i\leq k$ and $X_{k+1}=V(T_0)$.
Then $\X$ is an admissible 1-thin cover of $G-L'$ whose looped member is the empty set.
We have $f(T_0)=2|X_{k+1}|-3$, $f(T_i)=2|X_i|-3$ for all $1\leq i\leq k$ and similar to the previous case we can deduce

$$r^{lc}_2(G)=|L'|+\sum_{i=0}^kf(T_i)=|L'|+2|\emptyset|+\sum_{i=1}^{k+1}(2|X_i|-3)=\val(\X)$$
and this completes the proof.
\end{proof}
\begin{example}\label{exmp:rank_by_cover}
Let $G=(V,E,L)$ be the looped simple graph drawn in Figure \ref{fig:rank_example}. Let $X_0=\{v_1,v_2,v_3,v_7\}$, $X_1=\{v_3,v_4,v_5,v_6\}$, $X_2=\{v_8,v_9\}$ and $L'=\{l_1,l_2\}$.
Then the family $\X=\{X_0,X_1,X_2\}$ is an admissible 1-thin cover of $G-L'$ with $X_0$ being the looped member of $\X$ and Theorem \ref{thm:rank_vertex} gives
$$r^{lc}_2(G)\leq \val(\X)=|L'|+2|X_0|+\sum_{i=1}^2(2|X_i|-3)=2+8+5+1=16.$$
One can observe that we indeed have $r^{lc}_2(G)=16$ by calculating the rank of the rigidity matrix $R(G,p,q)$ of a generic linearly constrained framework $(G,p,q)$ in $\R^2$.
\end{example}
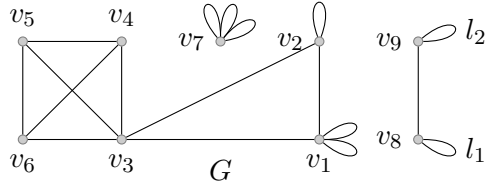
\begin{figure}[h]
\begin{center}
\begin{tikzpicture}[scale=1.3]
\node[roundnode] at (1,0) (a) [label=below:$v_1$]{}
	edge[in=0,out=320,loop]()
	edge[in=0,out=40,loop]();
\node[roundnode] at (-1,0) (b) [label=below:$v_3$]{}
	edge[] (a) ;
\node[roundnode] at (1,1) (c) [label=left:$v_2$]{}
	edge[] (a)
	edge[](b)
	edge[in=70,out=110,loop] ();
\node[roundnode] at (0,1) (c) [label=left:$v_7$]{}
	edge[in=90,out=130,loop]()
	edge[in=50,out=90,loop]()
	edge[in=10,out=50,loop]();
\node[roundnode] at (-1,1) (s) [label=above:$v_4$]{}
	edge[] (b);
\node[roundnode] at (-2,1) (t) [label=above:$v_5$]{}
	edge[] (b)
	edge[] (s);
\node[roundnode] at (-2,0) (z) [label=below:$v_6$]{}
	edge[] (b)
	edge[] (s)
	edge[] (t);
\node[roundnode] at (2,0) (v8) [label=left:$v_8$]{}
	edge[in=0,out=320,loop] node[pos=.6,right]{$l_1$}();
\node[roundnode] at (2,1) (v9) [label=left:$v_9$]{}
	edge[](v8)
	edge[in=0,out=40,loop] node[pos=.6,right]{$l_2$} ();
\node[] at (0,0) (G) [label=below:$G$] {};
\end{tikzpicture}
\end{center}
\caption{The looped simple graph $G$ in Example \ref{exmp:rank_by_cover}.}
\label{fig:rank_example}
\end{figure}

The proof of Case 2 above gives the relation between the covers we use and the covers used by Lov\'asz and Yemini.
Basically, if we would like to get the rank of a looped simple graph $G=(V,E,L)$, for which $L=\emptyset$ (a simple graph),
we should set the looped member as the empty set. Such covers are the 1-thin covers Lov\'asz and Yemini used to characterise
the rank function $r_2$ of $\calR_2$.

\section{$k$-balanced graphs and Rigidity}
\begin{figure}[ht]
\begin{center}

\begin{tikzpicture}[font=\small]

\node[roundnode] at (0cm,0cm) (u) [] {};
\node[roundnode] at (1cm,0cm) (v) [] {}
	edge[in=250,out=290,loop]()
	edge[](u);
\node[roundnode] at (1cm,1cm) (w) [] {}
	edge[in=70,out=110,loop]()
	edge[](u)
	edge[](v);
\node[roundnode] at (0cm,1cm) (z) [] {}
	edge[in=70,out=110,loop]()
	edge[](u)
	edge[](v)
	edge[](w);

\begin{scope}[xshift=4cm]
\node[roundnode] at (0cm,0cm) (u) [] {}
	edge[in=250,out=290,loop]();
\node[roundnode] at (1cm,0cm) (v) [] {}
	edge[in=250,out=290,loop]()
	edge[](u);
\node[roundnode] at (1cm,1cm) (w) [] {}
	edge[in=70,out=110,loop]()
	edge[](v);
\node[roundnode] at (0cm,1cm) (z) [] {}
	edge[in=70,out=110,loop]()
	edge[](u)
	edge[](w);
\end{scope}

\begin{scope}[xshift=8cm]
\node[roundnode] at (0cm,0cm) (u) [] {}
	edge[in=250,out=290,loop]();
\node[roundnode] at (1cm,0cm) (v) [] {}
	edge[in=250,out=290,loop]()
	edge[](u);
\node[roundnode] at (1cm,1cm) (w) [] {}
	edge[in=70,out=110,loop]()
	edge[](v);
\node[roundnode] at (0cm,1cm) (z) [] {}
	edge[in=70,out=110,loop]()
	edge[](u)
	edge[](w);
\begin{scope}[xshift=2.2cm]
\node[roundnode] at (0cm,0cm) (x) [] {}
	edge[in=250,out=290,loop]();
\node[roundnode] at (1cm,0cm) (y) [] {}
	edge[in=250,out=290,loop]()
	edge[](x);
\node[roundnode] at (1cm,1cm) (t) [] {}
	edge[in=70,out=110,loop]()
	edge[](y);
\node[roundnode] at (0cm,1cm) (s) [] {}
	edge[in=70,out=110,loop]()
	edge[](x)
	edge[](t);
\end{scope}
\end{scope}

\end{tikzpicture}

\end{center}
\caption{Some examples of 3-balanced graphs.}
\label{fig:looped_connectivity_examples}
\end{figure}
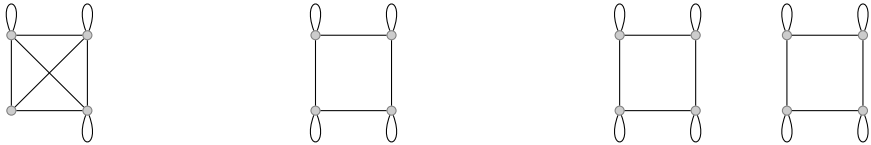

Let $G=(V,E,L)$ be a looped simple graph. We say $G$ is {\em $k$-balanced} if every connected component of $G-T$ where $T\subseteq V$ with $|T|\leq k$,
has at least $k-|T|$ vertices with loops. Note that by definition every $k$-balanced graph has at least
$k$ loops (by taking $T=\emptyset$). Note also that a $k$-balanced graph does not need to be connected. See Figure \ref{fig:looped_connectivity_examples}
for some examples of 3-balanced graphs. The underlying simple graph of the graph drawn on the left in Figure \ref{fig:looped_connectivity_examples} is 3-connected,
and adding three loops to different vertices makes it 3-balanced. (Adding $k$ loops to distinct vertices in a $k$-connected graph always gives a $k$-balanced graph.)
The graph in the middle is not 3-connected, but it is 3-balanced.
If we take the disjoint union of two copies of the graph in the middle, we obtain the graph on the right. Clearly, the graph we obtain after this operation is not even connected
whereas being 3-balanced is preserved.

The theorem below gives the relation between $k$-balancedness and rigidity of a linearly constrained framework in $\R^2$.
\begin{thm}\label{thm:minus3edges}
Let $G=(V,E,L)$ be a 6-balanced looped simple graph and $F\subset E\cup L$ with $|F|\leq 3$.
Then $G-F$ is rigid as a linearly constrained framework in $\R^2$.
\end{thm}
\begin{proof}
Let $G-F$ be a counterexample with $|V|+|F|+|L|$ being minimum and with respect to this $|E|$ being maximum. Let $F_e=F\cap E$, $F_l=F\cap L$ and
$H=G-F=(V,E\setminus F_e,L\setminus F_l)$.
Combining the fact that $H=G-F$ is not rigid and Theorem \ref{thm:rank_vertex}, there exists a set $L'\subseteq (L\setminus F_l)$
and an admissible 1-thin cover $\X=\{X_0,X_1,\ldots,X_k\}$ of $H-L'$ such that 
$$\val(\X)=|L'|+2|X_0|+\sum_{i=1}^k(2|X_i|-3)<2|V|.$$
Let $e_i=x_iy_i$, $1\leq i\leq |F_e|$ denote the edges in $F_e$ and put $X_{k+i}=\{x_i,y_i\}$, for all $1\leq i\leq |F_e|$.
The minimality of $|F|$ implies $X_{k+i}\not\subseteq X_j$ for all $1\leq i\leq |F_e|$ and for all $1\leq j\leq k$.
It also implies that the loops in $F_l$ are not incident with the vertices in $X_0$. Then the family $\X'=\X\cup\{X_{k+i}:1\leq i\leq |F_e|\}$
is an admissible 1-thin cover of $G-(L'\cup F_l)$ whose value satisfies
\begin{align}\label{eq:rank}
\begin{split}
\val(\X')&=|L'\cup F_l|+2|X_0|+\sum_{i=1}^{k+|F_e|}(2|X_i|-3)\\
&=\val(\X)+|F_l|+|F_e|=\val(\X)+|F|\leq\val(\X)+3<2|V|+3.
\end{split}
\end{align}
Note that the maximality of $|E|$ implies $(X_i,E_G(X_i))$ is a complete graph for all $0\leq i\leq k+|F_e|$ and the minimality of
$|L|$ implies every vertex $v\in V$ has at most one loop in $G$.
\begin{claim}\label{cla:two_Xi}
Every vertex $v\in V$ without a loop is contained in at least two $X_i$, $0\leq i\leq k+|F_e|$.
\end{claim}
\begin{proof}[Proof of Claim.]
Suppose the contrary and let $v$ be a loopless vertex that belongs to only one $X_j$ for some $0\leq j\leq k+|F_e|$.
The fact that $v$ has no loops and $G$ being 6-balanced imply that $d_G(v)\geq 6$. Since $X_j$ is the only set in the cover
$\X'$ that contains $v$, all neighbours of $v$ belong to $X_j$ and this implies $|X_j|\geq 7$. Therefore $1\leq j\leq k$ as $|X_{k+i}|=2$ for all $1\leq i\leq |F_e|$, and
so $v$ is not incident with the edges or loops in $F$.

Consider the graphs $\bar{G}=G-v$ and $\bar H=\bar G-F$, and let $\bar\X=\{\bar X_0,\bar X_1,\ldots,\bar X_k\}$ where $\bar X_j=X_j-v$ and $\bar X_i=X_i$ for all $i\neq j$.
Then $\bar\X$ is an admissible 1-thin cover of $\bar H-L'$ and we have
\begin{align*}
\val(\bar\X)=|L'|+2|\bar X_0|+\sum_{i=1}^k{(2|\bar X_i|-3)}&=|L'|+2|X_0|+\sum_{i=1}^k{(2|X_i|-3)}-2\\
									   &<2|V|-2=2|V-v|.
\end{align*}
Thus $\bar H=\bar G-F$ is not rigid as a linearly constrained framework in $\R^2$ by Theorem \ref{thm:rank_vertex}. By the minimality of $|V|$, the graph
$\bar G$ cannot be 6-balanced. Thefore there exists a non-empty set $T\subset (V-v)$ with $|T|\leq 5$ such that some
connected component(s) of $\bar G-T$ has less than $6-|T|$ loops. Let $C$ denote such a component. The facts that all neighbours of $v$ are contained in $X_j$ and
$(X_j,E_G(X_j))$ is a complete graph imply that all neighbours of $v$ belong to the same component $\bar C$ in $\bar G-T$. This implies
$C$ or $G[V(C)\cup\{v\}]$ is a connected component of $G-T$ with less than $6-|T|$ loops when $C\neq \bar C$ or $C=\bar C$, respectively and this gives a contradiction.
\end{proof}
Let $Y$ denote the set of vertices in $X_0$ whose neighbourhood is contained in $X_0$, and put $Z=X_0\setminus Y$ (i.e., $Z$ is the set of vertices in $X_0$ with at
least one neighbour in $V\setminus X_0$). For a vertex $v\in V$ let $l(v)$ denote the number of loops in $G$ incident with $v$. Since every vertex in $V$ has
at most one loop by the minimality of $|L|$, we have $0\leq l(v)\leq 1$ and $\sum_{v\in V\setminus X_0}l(v)=|L'\cup F_l|$.

\begin{claim}\label{cla:two_diff}
For $v\in V$ and $X_i\in \X'$, $0\leq i\leq k+|F_e|$, we have
\begin{itemize}
\item[](i) $\displaystyle{\sum_{X_i:v\in X_i}(2-\frac{3}{|X_i|})\geq 2-\frac{l(v)}{2}}$ when $v\in (V\setminus X_0)$,
\item[](ii) $\displaystyle{\sum_{X_i:v\in X_i}(2-\frac{3}{|X_i|})=2-\frac{3}{|X_0|}}$ when $v\in Y$, and
\item[](iii) $\displaystyle{\sum_{X_i:v\in X_i}(2-\frac{3}{|X_i|})\geq 2-\frac{3}{|X_0|}+\frac{1}{2}}$ when $v\in Z$.
\end{itemize}
\end{claim}
\begin{proof}[Proof of Claim]
\textbf{(i):} By relabelling if necessary, we may assume $v$ is contained in the sets $X_1,X_2,\ldots,X_m$ such that $|X_1|\geq |X_2|\geq\cdots\geq |X_m|$.

We first consider the case $l(v)=0$. Then by Claim \ref{cla:two_Xi}, $m\geq 2$. The facts that $G$ is 6-balanced and that $l(v)=0$ imply $d(v)\geq 6$.
Combining $|N(v)|=d(v)\geq 6$ with the fact that every edge incident with $v$ is contained in some $X_i$, we obtain $\sum_{X_i:v\in X_i}(|X_i|-1)\geq 6$.
Then depending on $m$, there are three subcases.\\
$m\geq 4$: Since $|X_i|\geq 2$, $\displaystyle{\sum_{X_i:v\in X_i}(2-\frac{3}{|X_i|})\geq 4\cdot \frac{1}{2}}=2$.\\
$m=3$: $\displaystyle{\sum_{X_i:v\in X_i}(|X_i|-1)\geq 6}$ implies $|X_1|\geq 3$. Then $\displaystyle{\sum_{X_i:v\in X_i}(2-\frac{3}{|X_i|})\geq 1+\frac{1}{2}+\frac{1}{2}=2}$.\\
$m=2$: $\displaystyle{\sum_{X_i:v\in X_i}(|X_i|-1)\geq 6}$ implies $|X_1|\geq 4$. Then either $|X_1|=|X_2|=4$ or $|X_1|=5$, $|X_2|\geq 3$ or $|X_1|\geq 6$, $|X_2|\geq 2$.
Then $\displaystyle{\sum_{X_i:v\in X_i}(2-\frac{3}{|X_i|})}$ is at least $\frac{5}{4}+\frac{5}{4}$ or $\frac{7}{5}+1$ or $\frac{3}{2}+\frac{1}{2}$, respectively.
Hence $\displaystyle{\sum_{X_i:v\in X_i}(2-\frac{3}{|X_i|})\geq 2-\frac{l(v)}{2}=2}$ is satisfied.

Now let us consider the case $l(v)=1$. Then $|N(v)|\geq 5$ since $G$ is 6-balanced. Similarly as above we can obtain $\sum_{X_i:v\in X_i}(|X_i|-1)\geq 5$.
Then there are three subcases depending on $m$.\\
$m\geq 3$: Since $|X_i|\geq 2$, $\displaystyle{\sum_{X_i:v\in X_i}(2-\frac{3}{|X_i|})\geq 3\cdot\frac{1}{2}= \frac{3}{2}}$.\\
$m=2$: $\displaystyle{\sum_{X_i:v\in X_i}(|X_i|-1)\geq 5}$ implies $|X_1|\geq 4$. Then $\displaystyle{\sum_{X_i:v\in X_i}(2-\frac{3}{|X_i|})\geq \frac{5}{4}+\frac{1}{2}\geq \frac{3}{2}}$.\\
$m=1$: Then $|X_1|\geq 6$ and $\displaystyle{\sum_{X_i:v\in X_i}(2-\frac{3}{|X_i|})=2-\frac{3}{|X_1|}\geq\frac{3}{2}}$.\\
Therefore $\displaystyle{\sum_{X_i:v\in X_i}(2-\frac{3}{|X_i|})\geq 2-\frac{l(v)}{2}=\frac{3}{2}}$ as required.\\
\textbf{(ii):} Follows from the fact that $X_0$ is the only member of $\X'$ which contains $v$. (This is implied by the facts that $v\in Y\subseteq X_0$, that is,
$(N(v)\cup\{v\})\subseteq X_0$, $(X_i,E_G(X_i))$ is complete for all $0\leq i\leq k+|F_e|$ and $\X'$ is 1-thin.)\\
\textbf{(iii):} Since $v\in Z\subset X_0$, there is an edge $vx$ for some $x\in V\setminus X_0$. Since $\X'$ covers the edge $vx$, there exists an $X_i$ with $x,v\in X_i$, for some
$1\leq i\leq k+|F_e|$. Hence $v$ is contained in at least two members of $\X'$ such that one of these members is $X_0$. Then the statement follows from the fact that $|X_i|\geq 2$
for all $1\leq i\leq k+|F_e|$.
\end{proof}
Note that since the loops in $L'\cup F_l$ are incident with vertices in $V\setminus X_0$ and $G$ is 6-balanced, $|Z|+|L'\cup F_l|\geq 6$ holds.
We now split the proof into two cases depending on $k+|F_e|$.\\
\textbf{Case 1.} $k+|F_e|=0$ and so $\X'=\{X_0\}$.\\
Then we have $V=X_0$. Since $(X_0,E_G(X_0))$ is a complete graph and $G$ has at least six vertices with loops we see that $G-F$ is rigid as a linearly constrained
framework in $\R^2$, a contradiction.\\
\textbf{Case 2.} $k+|F_e|\geq 1$.\\
In the following two subcases we will show $\val(\X')\geq 2|V|+3$. This together with (\ref{eq:rank}) will give us a contradiction and complete the proof.\\
\textbf{Case 2.1.} $X_0= \emptyset$.\\
Then we have
\begin{align*}
\val(\X')&=|L'\cup F_l|+\sum_{i=1}^{k+|F_e|}(2|X_i|-3)=|L'\cup F_l|+\sum_{i=1}^{k+|F_e|}|X_i|(2-\frac{3}{|X_i|})\\
								  &=|L'\cup F_l|+\sum_{v\in V}\sum_{X_i:v\in X_i}(2-\frac{3}{|X_i|})\\
								  &\geq |L'\cup F_l|+2|V|-\frac{|L'\cup F_l|}{2}\\
								  &=2|V|+\frac{|L'\cup F_l|}{2}\geq 2|V|+3,
\end{align*}
where the first inequality follows from Claim \ref{cla:two_diff}(i) and the last inequality follows from the facts that $|Z|+|L'\cup F_l|\geq 6$ and $Z\subset X_0=\emptyset$.\\
\textbf{Case 2.2.} $X_0\neq \emptyset$.\\
Then we have
\begin{align*}
\val(\X')&=|L'\cup F_l|+2|X_0|+\sum_{i=1}^{k+|F_e|}(2|X_i|-3)=|L'\cup F_l|+\sum_{i=0}^{k+|F_e|}(2|X_i|-3)+3\\
								  &=|L'\cup F_l|+\sum_{i=0}^{k+|F_e|}|X_i|(2-\frac{3}{|X_i|})+3\\
								  &=|L'\cup F_l|+\sum_{v\in V}\sum_{X_i:v\in X_i}(2-\frac{3}{|X_i|})+3\\
								  &=|L'\cup F_l|+\sum_{v\in Y}\sum_{X_i:v\in X_i}(2-\frac{3}{|X_i|})+\sum_{v\in Z}\sum_{X_i:v\in X_i}(2-\frac{3}{|X_i|})+\sum_{v\in V\setminus X_0}\sum_{X_i:v\in X_i}(2-\frac{3}{|X_i|})+3\\
								  &\geq |L'\cup F_l|+ 2|Y|-\frac{3|Y|}{|X_0|}+ 2|Z|-\frac{3|Z|}{|X_0|}+\frac{|Z|}{2}+2|V\setminus X_0|-\frac{|L'\cup F_l|}{2}+3\\
								  &=2|V|-\frac{3(|Y|+|Z|)}{|X_0|}+\frac{|Z|}{2}+\frac{|L'\cup F_l|}{2}+3=2|V|+\frac{|Z|+|L'\cup F_l|}{2}\geq 2|V|+3
\end{align*}
where the first inequality follows from Claim \ref{cla:two_diff}, and the last inequality follows from the fact that $|Z|+|L'\cup F_l|\geq 6$.
\end{proof}

Before giving our final result we need to give some definitions and some previous results first.
Let $G=(V,E,L)$ be a looped simple graph. We say that $G$ is {\em redundantly rigid} as a linearly contrained framework in $\R^2$
if $G-f$ is rigid as a linearly constrained framework in $\R^2$ for all $f\in E\cup L$.
Two linearly constrained frameworks $(G,p,q)$ and $(G,\bar p,q)$ are called {\em equivalent} if\\
$\bullet$ $||p(u)-p(v)||^2=||\bar p(u)-\bar p(v)||^2$ for all $u,v\in V$, and\\
$\bullet$ $p(v)\cdot q(l)=\bar p(v)\cdot q(l)$ for all incident pairs $v\in V$ and $l\in L$.

We say that $(G,p,q)$ is {\em globally rigid} if it has no equivalent realisation other than itself. In \cite{GJN}, it is shown that
global rigidity is a generic property and we have the following characterisation of global rigidity of linearly constrained frameworks in $\R^2$.
\begin{thm}\label{thm:glob_char}\cite{GJN}
Suppose $(G,p,q)$ is a generic linearly constrained framework in $\R^2$. Then $(G,p,q)$ is globally rigid if and only if
\begin{itemize}
\item[](i) each connected component of $G$ is either a single vertex with two loops or is redundantly rigid, and
\item[](ii) each connected component of $G-X$ has at least one loop for all $X\subseteq V(G)$ with $|X|=2$.
\end{itemize}
\end{thm}
\noindent \textbf{Remark.} The graphs that satisfy (ii) in the statement of Theorem \ref{thm:glob_char} are called
``2-balanced'' (or ``balanced'' for short) in \cite{GJN}. Since this terminology is slightly different from the definition of being ``$k$-balanced'' in this paper,
we wrote this condition explicitly in order not to cause a confusion.

Let $G=(V,E,L)$ be a 6-balanced looped simple graph. By Theorem \ref{thm:minus3edges}, $G$ is redundantly rigid as a linearly constrained framework in $\R^2$.
Moreover, it is easy to see that as being 6-balanced $G$ satisfies (ii) in the statement of Theorem \ref{thm:glob_char}.
Therefore we immediately obtain the following result.
\begin{thm}
Every generic realisation of a 6-balanced looped simple graph as a linearly constrained framework in $\R^2$ is globally rigid.
\end{thm}
\section{Further Remarks and Examples}

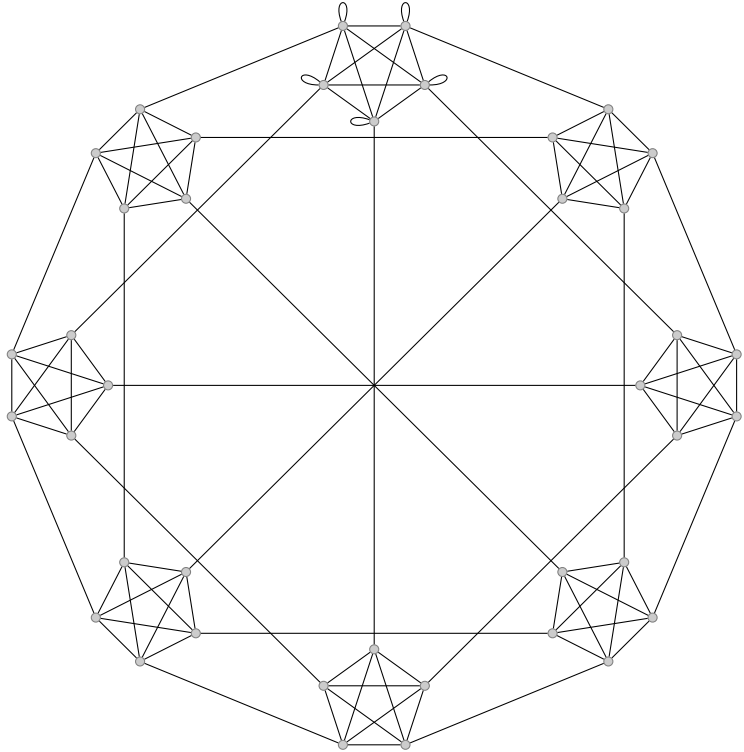
\begin{figure}[ht]
\begin{center}

\begin{tikzpicture}[font=\small,scale=0.7]

\begin{scope}[rotate=0]
\begin{scope}[xshift=6cm]
\begin{scope}[rotate=-90]
\node[roundnode] at (-18:1cm) (a1) [] {};
\node[roundnode] at (54:1cm) (a2) [] {}
	edge[] (a1);
\node[roundnode] at (126:1cm) (a3) [] {}
	edge[] (a1)
	edge[] (a2);
\node[roundnode] at (198:1cm) (a4) [] {}
	edge[] (a1)
	edge[] (a2)
	edge[] (a3);
\node[roundnode] at (270:1cm) (a5) [] {}
	edge[] (a1)
	edge[] (a2)
	edge[] (a3)
	edge[] (a4);
\end{scope}
\end{scope}
\end{scope}
\begin{scope}[rotate=45]
\begin{scope}[xshift=6cm]
\begin{scope}[rotate=-90]
\node[roundnode] at (-18:1cm) (b1) [] {};
\node[roundnode] at (54:1cm) (b2) [] {}
	edge[] (b1)
	edge[] (a3);
\node[roundnode] at (126:1cm) (b3) [] {}
	edge[] (b1)
	edge[] (b2);
\node[roundnode] at (198:1cm) (b4) [] {}
	edge[] (b1)
	edge[] (b2)
	edge[] (b3);
\node[roundnode] at (270:1cm) (b5) [] {}
	edge[] (b1)
	edge[] (b2)
	edge[] (b3)
	edge[] (b4);
\end{scope}
\end{scope}
\end{scope}

\begin{scope}[rotate=90]
\begin{scope}[xshift=6cm]
\begin{scope}[rotate=-90]
\node[roundnode] at (-18:1cm) (c1) [] {}
	edge[] (a4)
	edge[in=0,out=40,loop]();
\node[roundnode] at (54:1cm) (c2) [] {}
	edge[] (c1)
	edge[] (b3)
	edge[in=70,out=110,loop]();
\node[roundnode] at (126:1cm) (c3) [] {}
	edge[] (c1)
	edge[] (c2)
	edge[in=70,out=110,loop]();
\node[roundnode] at (198:1cm) (c4) [] {}
	edge[] (c1)
	edge[] (c2)
	edge[] (c3)
	edge[in=180,out=140,loop]();
\node[roundnode] at (270:1cm) (c5) [] {}
	edge[] (c1)
	edge[] (c2)
	edge[] (c3)
	edge[] (c4)
	edge[in=200,out=160,loop]();
\end{scope}
\end{scope}
\end{scope}

\begin{scope}[rotate=135]
\begin{scope}[xshift=6cm]
\begin{scope}[rotate=-90]
\node[roundnode] at (-18:1cm) (d1) [] {}
	edge[] (b4);
\node[roundnode] at (54:1cm) (d2) [] {}
	edge[] (d1)
	edge[] (c3);
\node[roundnode] at (126:1cm) (d3) [] {}
	edge[] (d1)
	edge[] (d2);
\node[roundnode] at (198:1cm) (d4) [] {}
	edge[] (d1)
	edge[] (d2)
	edge[] (d3);
\node[roundnode] at (270:1cm) (d5) [] {}
	edge[] (d1)
	edge[] (d2)
	edge[] (d3)
	edge[] (d4);
\end{scope}
\end{scope}
\end{scope}

\begin{scope}[rotate=180]
\begin{scope}[xshift=6cm]
\begin{scope}[rotate=-90]
\node[roundnode] at (-18:1cm) (e1) [] {}
	edge[] (c4);
\node[roundnode] at (54:1cm) (e2) [] {}
	edge[] (e1)
	edge[] (d3);
\node[roundnode] at (126:1cm) (e3) [] {}
	edge[] (e1)
	edge[] (e2);
\node[roundnode] at (198:1cm) (e4) [] {}
	edge[] (e1)
	edge[] (e2)
	edge[] (e3);
\node[roundnode] at (270:1cm) (e5) [] {}
	edge[] (e1)
	edge[] (e2)
	edge[] (e3)
	edge[] (e4)
	edge[] (a5);
\end{scope}
\end{scope}
\end{scope}

\begin{scope}[rotate=225]
\begin{scope}[xshift=6cm]
\begin{scope}[rotate=-90]
\node[roundnode] at (-18:1cm) (f1) [] {}
	edge[] (d4);
\node[roundnode] at (54:1cm) (f2) [] {}
	edge[] (f1)
	edge[] (e3);
\node[roundnode] at (126:1cm) (f3) [] {}
	edge[] (f1)
	edge[] (f2);
\node[roundnode] at (198:1cm) (f4) [] {}
	edge[] (f1)
	edge[] (f2)
	edge[] (f3);
\node[roundnode] at (270:1cm) (f5) [] {}
	edge[] (f1)
	edge[] (f2)
	edge[] (f3)
	edge[] (f4)
	edge[] (b5);
\end{scope}
\end{scope}
\end{scope}

\begin{scope}[rotate=270]
\begin{scope}[xshift=6cm]
\begin{scope}[rotate=-90]
\node[roundnode] at (-18:1cm) (g1) [] {}
	edge[] (e4);
\node[roundnode] at (54:1cm) (g2) [] {}
	edge[] (g1)
	edge[] (f3);
\node[roundnode] at (126:1cm) (g3) [] {}
	edge[] (g1)
	edge[] (g2);
\node[roundnode] at (198:1cm) (g4) [] {}
	edge[] (g1)
	edge[] (g2)
	edge[] (g3)
	edge[] (a1);
\node[roundnode] at (270:1cm) (g5) [] {}
	edge[] (g1)
	edge[] (g2)
	edge[] (g3)
	edge[] (g4)
	edge[] (c5);
\end{scope}
\end{scope}
\end{scope}

\begin{scope}[rotate=315]
\begin{scope}[xshift=6cm]
\begin{scope}[rotate=-90]
\node[roundnode] at (-18:1cm) (h1) [] {}
	edge[] (f4);
\node[roundnode] at (54:1cm) (h2) [] {}
	edge[] (h1)
	edge[] (g3);
\node[roundnode] at (126:1cm) (h3) [] {}
	edge[] (h1)
	edge[] (h2)
	edge[] (a2);
\node[roundnode] at (198:1cm) (h4) [] {}
	edge[] (h1)
	edge[] (h2)
	edge[] (h3)
	edge[] (b1);
\node[roundnode] at (270:1cm) (h5) [] {}
	edge[] (h1)
	edge[] (h2)
	edge[] (h3)
	edge[] (h4)
	edge[] (d5);
\end{scope}
\end{scope}
\end{scope}
\end{tikzpicture}

\end{center}
\caption{A non-rigid 5-balanced looped simple graph.}
\label{fig:5-balanced}
\end{figure}

We cannot replace the number 6 (i.e.,\ being 6-balanced) in Theorem \ref{thm:minus3edges} by a smaller number. To see this consider the
looped simple graph $G=(V,E,L)$ in Figure \ref{fig:5-balanced}. Let $H$ be the underlying simple graph of $G$.
Since $H$ is 5-connected and adding a single loop to five distinct vertices turns $H$ into $G$, the graph $G$ is 5-balanced.
Lov\'asz and Yemini \cite{LY} used Theorem \ref{thm:LY} to show that the graph $H$ is not  rigid as a bar-and-joint framework in $\R^2$. 
Using Theorem \ref{thm:rank_vertex} in a similar way, we can show that the graph $G$ is not rigid
as a linearly constrained framework in $\R^2$. To see this let $X_0$ be the set of vertices with loops, let $X_1,X_2,\ldots ,X_7$
be the vertex sets of copies of the other $K_5$'s, and let $X_8,X_9,\ldots,X_{27}$ be the sets of the endpoints of the simple edges
that connect distinct copies of $K_5$'s. Then $\X=\{X_0,X_1,\ldots,X_{27}\}$ is an admissible 1-thin cover of $G$ ($L'=\emptyset$).
Thus Theorem \ref{thm:rank_vertex} gives 
$$r^{lc}_2(G)\leq \val(\X)=|L'|+2|X_0|+\sum_{i=1}^{27}(2|X_i|-3)=0+10+7\cdot 7+20= 79<80=2|V|$$
implying that $G$ is not rigid as a linearly constrained framework in $\R^2$.

If we remove more than three edges or loops from a 6-balanced looped simple graph, we may end up with a non-rigid graph.
To see this let $G=(V,E,L)$ be a graph obtained from a 6-connected simple graph by adding a single loop to six distinct vertices. Then clearly, $G$ is 6-balanced.
Let $l_1,l_2\ldots,l_6$ denote the loops of $G$. Consider the graph $H=G-F$ where $F=\{l_1,l_2,l_3,l_4\}$, and let $L'=\{l_5,l_6\}$, $X_0=\emptyset$ and $X_1=V$.
Then $\X=\{X_0,X_1\}$ is an admissible 1-thin cover of $H-L'$ whose looped member is the empty set. Thus Theorem \ref{thm:rank_vertex} gives
$$r^{lc}_2(H)\leq \val(\X)=|L'|+2|X_0|+2|X_1|-3=2+0+2|V|-3=2|V|-1<2|V|$$
implying that $H$ is not rigid as a linearly constrained framework in $\R^2$.

\noindent\textbf{Acknowledgments.} We would like to thank Bill Jackson for his useful comments on an earlier version of the manuscript.
We would also like to thank Anthony Nixon for reminding the relation between 6-balancedness and global rigidity, and Matteo Gallet for some corrections.
Lastly, we would like to thank the anonymous referees whose comments helped to improve this paper.


\begin{thebibliography}{99}

\bibitem{A}
{T.G.\ Abbot},  {Generalizations of Kempe's Universality Theorem}.
MSc thesis, MIT (2008).
http://web.mit.edu/tabbott/www/papers/mthesis.pdf

\bibitem{AR}
{L.\ Asimow and B.\ Roth}, { The rigidity of graphs}, { Trans.~Am.~Math.~Soc.}
{245} (1978), 279-289.

\bibitem{CGJN}
{J.\ Cruickshank, H.\ Guler, B.\ Jackson and A.\ Nixon},
{Rigidity of Linearly Constrained Frameworks}, {International Mathematics Research Notices}
{12} (2020), 3824–3840.

\bibitem{Edm}
{J.\ Edmonds}, {Submodular functions, matroids, and certain polyhedra}, {in: Combinatorial Structures and their Applications, R.\ Guy, H.\ Hanani, N.\ Sauer, and J.\
Schönheim (Eds.)}, Gordon and Breach, New York, (1970), pp. 69-87.

\bibitem{AF}
{A.\ Frank}, { Connections in Combinatorial Optimization}, {Oxford Lecture Series in Mathematics and Its Applications 38}, {Oxford University Press}, New York, 2011.

\bibitem{GJN}
{H.\ Guler, B.\ Jackson and A.\ Nixon},
{Global Rigidity of 2D Linearly Constrained Frameworks}, {International Mathematics Research Notices}
{22} (2021), 16811–16858.

\bibitem{JNT}
{B.\ Jackson, A.\ Nixon and S.-I.\ Tanigawa}, {An Improved Bound for the Rigidity of Linearly Constrained Frameworks},
{SIAM J. Discrete Math} 35(2) (2021) 928–933.

\bibitem{KaTa}
{N.\ Katoh and S.-I.\ Tanigawa}, {On the Infinitesimal Rigidity of Bar-and-Slider Frameworks}, {in: ISAAC 2009: Algorithms and Computation, Dong, Y., Du, DZ., Ibarra, O. (Eds.)}, {Lecture Notes in Computer Science, vol 5878.} Springer, Berlin, Heidelberg, (2009), pp. 524–533.

\bibitem{KatTan}
{N.\ Katoh and S.-I.\ Tanigawa}, {Rooted-tree decompositions with matroid constraints and the infinitesimal rigidity of frameworks with boundaries},
{SIAM J. Discrete Math.} 27(1) (2013), 155-185.

\bibitem{L} {G.\ Laman},
{On graphs and rigidity of plane skeletal structures},
{J. Engineering Math.} {4} (1970), 331-340

\bibitem{LY} {L.\ Lov\'asz and  Y.\ Yemini},
{On generic rigidity in the plane}, {SIAM J. Algebraic Discrete
Methods} {3}  (1982), 91-98.


\bibitem{PG} {H.\ Pollaczek-Geiringer}, {\"{U}ber die Gliederung ebener Fachwerke}, {Zeitschrift f\"{u}r Angewandte Mathematik und Mechanik (ZAMM)} {7}  (1927), 58-72.

\bibitem{ST} {I.\ Streinu and L.\ Theran}, {Slider-pinning rigidity: a Maxwell-Laman-type theorem}, Discrete and
Computational Geometry, {44} (2010) 812–837.

\end{thebibliography}
\end{document}